\documentclass[11pt]{article}
\usepackage{amssymb,amsfonts,amsmath, psfrag,eepic,colordvi}
\newtheorem{theo}{Theorem}

\makeatletter \@addtoreset{equation}{section}
\@addtoreset{theo}{}\makeatother

\def\qed{\hfill \rule{4pt}{7pt}}
\def\pf{\noindent {\it Proof.} }

\begin{document}
\begin{center}
{\LARGE\bf
Weighted $2$-Motzkin Paths}\\[12pt]
William Y.C. Chen$^1$, Sherry H.F. Yan$^2$ and   Laura L.M. Yang$^3$ \\[6pt]
Center for Combinatorics, LPMC\\
Nankai University,
Tianjin 300071, P. R. China\\[5pt]
$^1$chen@nankai.edu.cn, $^2$huifangyan@eyou.com,
$^3$yanglm@hotmail.com
\end{center}

\vskip 6mm

\noindent {\bf Abstract. }   This paper is motivated by two
problems recently proposed by Coker on combinatorial identities
related to the Narayana polynomials and the Catalan numbers. We
find that a bijection of Chen, Deutsch and Elizalde can be used to
provide combinatorial interpretations of the identities of Coker
when it is applied to weighted plane trees. For the sake of
presentation of our combinatorial correspondences, we provide a
description of the bijection of Chen, Deutsch and Elizalde in a
slightly different manner in the form of a direct construction
from plane trees to $2$-Motzkin paths without the intermediate
step involving the Dyck paths.

\vskip 8mm



\noindent
 {\bf Suggested Running Title:} Weighted $2$-Motzkin Paths

  \noindent
 {\bf Corresponding Author:} William Y. C. Chen, Email:
 chen@nankai.edu.cn

\vskip 1cm

\section{Introduction}

The structure of $2$-Motzkin paths, introduced by Barcucci, Lungo,
Pergola and Pinzani \cite{blpp}, has proved to be highly efficient
in the study of plane trees, Dyck paths, Motzkin paths,
noncrossing partitions, RNA secondary structures,
Devenport-Schinzel sequences, and combinatorial identities, see
\cite{ds, Klazar, SW}. While it is most natural to establish a
correspondence between Dyck paths of length $2n$ and $2$-Motzkin
paths of length $n-1$, Deutsch and Shapiro came to the realization
that direct correspondences between plane trees  and $2$-Motzkin
paths can have many applications. Recently, Chen, Deutsch and
Elizalde \cite{cde} found  bijections between plane trees and
$2$-Motzkin paths  for the enumeration of plane trees by the
numbers of old and young leaves \cite{cde}. The main result of
this paper is to show that the bijection of Chen, Deutsch and
Elizalde, presented in a slightly different manner, can be applied
to weighted plane trees in order to give combinatorial
interpretations of two identities involving the Narayana numbers
and Catalan numbers due to Coker \cite{coker}. This leads to the
solutions of the two open problems left in the paper \cite{coker}.

Recall that a {\it $2$-Motzkin path} is a lattice path starting at
$(0,0)$ and ending at $(n,0)$ but never going below the $x$-axis,
with possible steps $(1,1)$, $(1,0)$ and $(1,-1)$, where the level
steps $(1,0)$ can be either of two kinds: straight and wavy. The
{\it length} of the path is defined to the number of its steps.
Deutsch and Shapiro \cite{ds} presented a bijection between plane
trees with $n$ edges and $2$-Motzkin paths of length $n-1$. So the
number of $2$-Motzkin paths of length $n-1$ equals the Catalan
number
\[ C_n = {1 \over n+1} \, {2n \choose n}.\]

Recently, Coker \cite{coker} established very interesting
combinatorial identities ((\ref{q1}) and (\ref{q2}) below) by
using generating functions and the Lagrange inversion formula
based the study of multiple Dyck paths. A multiple Dyck path is a
lattice path starting at $(0,0)$ and ending at $(2n, 0)$ with big
steps that can be regarded as segments of consecutive up steps or
consecutive down steps in an ordinary Dyck path. Note that the
notion of multiple Dyck path is formulated by Coker in different
coordinates.  The main ingredients in Coker's identities are the
Catalan number and the Narayana numbers:
$$
N(n,k)=\frac{1}{n}{n \choose k}{n \choose k-1},
$$
which counts the number of all plane trees with $n$ edges and $k$
leaves \cite{SW,stanley2}. It is sequence $A001263$ in
\cite{sloane}. Coker \cite{coker} left the following two open
problems:

\noindent {\bf Problem 1.} Find a bijective proof of the following
identity
\begin{equation}\label{q1}
\sum_{k=1}^n\frac{1}{n}{n \choose k}{n \choose
k-1}4^{n-k}=\sum_{k=0}^{\lfloor{(n-1)}/{2}\rfloor}C_k{n-1 \choose
2k}4^k5^{n-2k-1}.
\end{equation}
\noindent {\bf Problem 2.} Find a combinatorial explanation for
the following identity
\begin{equation}\label{q2}
\sum_{k=1}^n\frac{1}{n}{n \choose k}{n \choose
k-1}x^{2k}(1+x)^{2n-2k}=x^2\sum_{k=0}^{n-1}C_{k+1}{n-1 \choose
k}x^k(1+x)^{k},
\end{equation}
which is equation (6.2) in \cite{coker}. The above identity
(\ref{q1}) is a special case of the following identity:
\begin{equation}\label{q1x}
\sum_{k=0}^n \, {1\over n} \, {n\choose k} {n \choose k-1} t^{n-k}
= \sum_{k=0}^{\lfloor (n-1)/2\rfloor}\, C_k\, {n-1\choose 2k}\,t^k
(1+t)^{n-2k-1},
\end{equation}
where the left hand side of (\ref{q1x}) is the Narayana
polynomial, as denoted by ${\cal N}_n(t)$  in \cite{coker}. The
identity (\ref{q1x}) is the relation (4.4) in \cite{coker}, which
can be derived as an identity on the Narayana numbers and the
Catalan numbers due to Simion and Ullman \cite{SU}, see also
\cite{CDD}.
 Remarkably, (\ref{q1x}) has many
consequences as pointed by Coker \cite{coker}. For example, it
implies the classical identity of Touchard \cite{T}, and the
formula on the little Schr\"{o}der numbers in terms of the Catalan
numbers \cite{G}. The reason for the evaluation of ${\cal N}_n(t)$
at $t=4$ lies in the fact that ${\cal N}_n(4)$ equals the number
$d_n$ of multiple Dyck paths of length $2n$. The first few values
of $d_n$ for $n=0, 1,2,3,4, 5,6,7$ are \[ 1,\;1,\; 5,\; 29,\;
185,\; 1257,\; 8925,\; 65445,\]
 which is the sequence $A059231$ in \cite{sloane}. From
the interpretation of Narayana numbers in terms of Dyck paths of
length $2n$ and of $k$ peaks, it is not difficult to show that
$d_n={\cal N}_n(4)$. However,  the right hand side of (\ref{q1})
does not seems to be obvious, which is obtained by establishing a
functional equation and by using the Lagrange inversion formula.
The natural question as raised by Coker \cite{coker} is to find a
combinatorial interpretation of (\ref{q1}). Note that the
enumeration of multiple Dyck paths has also been studied
independently by Sulanke \cite{Sulanke} and Woan \cite{Woan}.

The relation (\ref{q2}) was established  from the enumeration of
multiple Dyck paths of length $2n$ with a given number of  steps.
Let $\lambda_{n,j}$ be the number of multiple Dyck paths of length
$2n$ and $j$ steps, ${\cal P}_n(x)$ be the polynomial
\[ {\cal P}_n(x)= \sum_{j=2}^{2n} \, \lambda_{n,j}\, x^j.\]
It was shown that
\begin{equation}
{\cal P}_n(x) = \sum_{k=1}^n \, {1 \over n} {n\choose k} {n
\choose k-1} x^{2k} (1+x)^{2n-2k},
\end{equation}
which can be restated as
\begin{equation}
{\cal P}_n(x) = x^{2n} {\cal N}_n((1+x^{-1})^2).
\end{equation}
Coker \cite{coker} discovered the connection between ${\cal
P}_n(x)$ and the polynomial ${\cal R}_n(x)$ introduced by Denise
and Simion \cite{denise-s} in their study of the number of
exterior pairs of Dyck paths of length $2n$. The polynomials
${\cal R}_n(x)$ have the following expansion:
\begin{equation}
{\cal R}_n(x) = \sum_{k=0}^{n-1} \, (-1)^k C_{k+1}{n-1 \choose
k}\, x^k (1-x)^k.
\end{equation}
It now becomes clear that the identity (\ref{q2}) can be rewritten
as
\begin{equation}\label{PR}
{\cal P}_n(x) =x^2 {\cal R}_n(-x).
\end{equation}

To give  combinatorial interpretations of both (\ref{q1}) and
(\ref{q2}), we apply the bijection of Chen, Deutsch and Elizalde
\cite{cde} to weighted plane trees to get weighted $2$-Motzkin
paths. Then we use weight-preserving operations on $2$-Motzkin
paths to derive the desired combinatorial identities. More
precisely, these weight-preserving operations are essentially the
reductions from weighted $2$-Motzkin paths to Dyck paths and
$2$-Motzkin paths. It would be interesting to find a direct
correspondence on Dyck paths which leads to a combinatorial
interpretation of (\ref{PR}).

\section{Weighted $2$-Motzkin Paths}

Let us review a bijection between plane trees and $2$-Motzkin
paths due to Chen, Deutsch and Elizalde \cite{cde}, which is
devised for the enumeration of plane trees with $n$ edges and a
fixed number of old leaves and a fixed number of young leaves.
Such a consideration of old and young leaves reflects to the four
types of steps of $2$-Motzkin paths. For the purpose of this
paper, we present a  slightly modified version of the nonrecursive
bijection in \cite{cde}. Our terminology is also somewhat
different.

 For a plane tree $T$, a vertex of $v$ is
called a leaf if it does
 not have any children. An internal vertex is a
vertex that has at least one child. An edge is denoted as a pair
$(u, v)$  of vertices such that $v$ is a child of $u$.   Let $u$
be an internal vertex, and $v_1$, $v_2$, $\ldots$, $v_k$ be the
children of $u$ listed from left to right. Then we call $v_k$ an
{\it exterior vertex} and $(u, v_k)$ an {\it exterior edge}. If
$k>1$, then the edges $(u, v_1)$, $(u, v_2)$ $\ldots$, and $(u,
v_{k-1})$ are called {\it interior edges} and $v_1, v_2, \ldots,
v_{k-1}$ are called {\it interior vertices}. An edge containing a
leaf vertex is called a {\it terminal edge}. Let $u$ be the root
of $T$, $(u, u_1)$ be the exterior edge of $u$,
    $(u_1, u_2)$ be the exterior edge of $u_1$, and so on, finally $(u_{k-1}, u_{k})$
    be the exterior edge of $u_{k-1}$  such that $u_k$ is a leaf.
    The exterior edge $(u_{k-1}, u_k)$ is called the {\it critical edge} of $T$.
    To summarize, the edges of a plane tree $T$ are classified  into five
     categories.
\begin{itemize}

    \item Non-terminal interior edges.

    \item Non-terminal exterior edges.

    \item Terminal interior edges.

    \item Terminal exterior  edges (which do not include the
    critical edge).

    \item The critical edge.
\end{itemize}

Note that the {\it critical edge} of $T$ is the last encountered
edge when we traverse the edges of $T$ in preorder. From the above
classification on the edges of a plane tree, it is easy to
describe the Chen-Deutsch-Elizalde bijection between plane trees
and $2$-Motzkin paths by the preorder traversal of the edges of
$T$. To be precise, let $u$ be the root of $T$,  $v_1, v_2,
\ldots, v_k$ be the children of $u$, and $T_1, T_2, \ldots, T_k$
be the subtrees of $T$ rooted at $v_1, v_2, \ldots, v_k$,
respectively. Then the preorder traversal of the edges of $T$,
denoted by $P(T)$, is a linear order of the edges of $T$
recursively defined by
\[ (u, v_1) \, P(T_1) \, (u, v_2) \, P(T_2) \, \cdots \, (u, v_k)
\, P(T_k).\]

\noindent {\bf The Bijection of Chen, Deutsch and Elizalde
\cite{cde}:}
 Let $T$ be any nonempty plane tree.
At each step of the traversal of the edges of $T$ in preorder,
\begin{description}
    \item (i) draw an up step for a non-terminal interior edge;
    \item (ii) draw a straight level step for a non-terminal exterior edge;
    \item (iii) draw a wavy level step for a terminal interior edge;
    \item (iv) draw a down step for a  terminal exterior edge;
    \item (v) do nothing for the critical edge.
\end{description}
\begin{figure}[h,t]
\begin{center}
\begin{picture}(400,60)
\setlength{\unitlength}{1mm} \linethickness{0.4pt}
\linethickness{0.4pt} \put(5,15){\line(1,1){10}}
\put(15,25){\circle*{1}} \put(15,25){\line(-1,-2){5}}
\put(5,15){\circle*{1}} \put(10,15){\circle*{1}}
\put(5,10){\line(0,1){5}} \put(5,10){\line(1,-2){2.5}}
\put(5,10){\line(-1,-2){2.5}} \put(5,10){\line(0,-1){5}}
\put(5,10){\circle*{1}} \put(7.5,5){\circle*{1}}
\put(2.5,0){\circle*{1}} \put(5,5){\circle*{1}}
\put(5,5){\line(-1,-2){2.5}} \put(5,5){\line(1,-2){2.5}}
\put(2.5,5){\circle*{1}} \put(7.5,0){\circle*{1}}
\put(15,25){\line(0,-1){10}} \put(15,15){\line(1,-2){2.5}}
\put(15,15){\line(-1,-2){2.5}} \put(15,15){\circle*{1}}
\put(17.5,10){\circle*{1}} \put(12.5,10){\circle*{1}}

\put(17.5,10){\line(-1,-2){2.5}} \put(17.5,10){\line(1,-2){2.5}}
\put(20,5){\line(0,-1){5}} \put(15,5){\circle*{1}}
\put(20,5){\circle*{1}} \put(20,0){\circle*{1}}
\put(15,25){\circle*{1}} \put(15,25){\line(1,-1){10}}
\put(25,15){\circle*{1}} \put(25,15){\line(2,-1){10}}
\put(35,10){\circle*{1}} \put(25,15){\line(1,-1){5}}
\put(30,10){\circle*{1}} \put(25,15){\line(0,-1){5}}
\put(25,10){\circle*{1}} \put(35,10){\line(0,-1){5}}
\put(35,5){\circle*{1}} \put(30,10){\line(0,-1){5}}
\put(30,5){\circle*{1}}
\linethickness{1pt} \put(45,10){\vector(1,0){6}}
\put(51,9){\vector(-1,0){6}}
\linethickness{0.4pt} \put(60,5){\circle*{1}}
\put(60,5){\line(1,1){4}}

\put(64,9){\circle*{1}} \put(64,9){\line(1,0){4}}

\put(68,9){\circle*{1}} \multiput(68,9)(1,0){4}{\line(1,1){0.5}}
\multiput(68.5,9.5)(1,0){4}{\line(1,-1){0.5}}

\put(72,9){\circle*{1}} \put(72,9){\line(1,1){4}}

\put(76,13){\circle*{1}} \multiput(76,13)(1,0){4}{\line(1,1){0.5}}
\multiput(76.5,13.5)(1,0){4}{\line(1,-1){0.5}}

\put(80,13){\circle*{1}} \put(80,13){\line(1,-1){4}}

\put(84,9){\circle*{1}} \put(84,9){\line(1,-1){4}}

\put(88,5){\circle*{1}} \multiput(88,5)(1,0){4}{\line(1,1){0.5}}
\multiput(88.5,5.5)(1,0){4}{\line(1,-1){0.5}}

\put(92,5){\circle*{1}} \put(92,5){\line(1,1){4}}

\put(96,9){\circle*{1}} \multiput(96,9)(1,0){4}{\line(1,1){0.5}}
\multiput(96.5,9.5)(1,0){4}{\line(1,-1){0.5}}

\put(100,9){\circle*{1}} \put(100,9){\line(1,0){4}}

\put(104,9){\circle*{1}} \multiput(104,9)(1,0){4}{\line(1,1){0.5}}
\multiput(104.5,9.5)(1,0){4}{\line(1,-1){0.5}}

\put(108,9){\circle*{1}} \put(108,9){\line(1,0){4}}

\put(112,9){\circle*{1}} \put(112,9){\line(1,-1){4}}

\put(116,5){\circle*{1}} \put(116,5){\line(1,0){4}}

\put(120,5){\circle*{1}} \multiput(120,5)(1,0){4}{\line(1,1){0.5}}
\multiput(120.5,5.5)(1,0){4}{\line(1,-1){0.5}}

\put(124,5){\circle*{1}} \put(124,5){\line(1,1){4}}

\put(128,9){\circle*{1}} \put(128,9){\line(1,-1){4}}

\put(132,5){\circle*{1}} \put(132,5){\line(1,0){4}}
\put(136,5){\circle*{1}}
\end{picture}
\end{center}
\caption{The Bijection $\Phi$} \label{Upsilon}
\end{figure}

It is easy to see that we have obtained a $2$-Moztkin path. More
precisely,  a plane tree with $n$ edges corresponds to a
$2$-Motzkin path of length $n-1$. The above bijection is  denoted
by $\Phi$. As a hint to why the above bijection works, one may
check that for any plane tree $T$,
\begin{center}
 $\#$ non-terminal interior edges $=$
 $\#$ terminal exterior edges.
\end{center}

We are now ready to assign weights to the edges of a plane tree
$T$ in order to obtain  combinatorial interpretations of the
identities (\ref{q1}) and (\ref{q2}). The weights of the edges of
a plane tree will translate into weights of steps  of the
corresponding $2$-Motzkin path. In fact, we will take slightly
different formulations of  (\ref{q1x}) and (\ref{q2}).

\begin{theo}
For $n \geq 1$, we have
\begin{equation}\label{identity1}
\sum_{k=1}^n\frac{1}{n}{n \choose k}{n \choose
k-1}x^{k-1}=\sum_{k=0}^{\lfloor{(n-1)}/{2}\rfloor}C_k{n-1 \choose
2k}x^k(1+x)^{n-2k-1}.
\end{equation}
\end{theo}

\pf Let $T$ be a plane tree with $n$ edges. We assign the weights
to the edges of $T$ by the following rule: All the terminal edges
except the critical edge are given the weight $x$ and all other
edges are given the weight $1$. The weight of $T$ is the product
of the weights of its edges.  Then the left hand side of
(\ref{identity1}) is the sum of the weights of all plane trees
with $n$ edges.

By the above bijection $\Phi$, the set of weighted plane trees
with $n$ edges is mapped to the set of $2$-Motzkin paths of length
$n-1$ in which all the down steps and wavy level steps are given
the weight $x$, and other steps are given the weight $1$. Consider
the weighted $2$-Motzkin paths of length $n-1$ that have $k$ up
steps and $k$ down steps. These $k$ up steps and $k$ down steps
form a Dyck path of length $2k$. The binomial coefficient ${n-1
\choose 2k}$ comes from the choices of the $2k$ positions for
these $k$ up steps and $k$ down steps.  The remaining $n-2k-1$
steps are either wavy level steps or straight level steps. Since
only a wavy level step carries the weight $x$, the total
contributions of the weights of $n-2k-1$ level steps amount to
$(1+x)^{n-2k-1}$. The $k$ up steps would contribute $x^k$.
Therefore, the right hand side of (\ref{identity1}) equals the
total contributions of all the weighted $2$-Motzkin paths of
length $n-1$, as desired. \qed

\begin{theo}
For $n \geq 1$, we have
\begin{equation}\label{identity2}
\sum_{k=1}^n\frac{1}{n}{n \choose k}{n \choose
k-1}x^{2(k-1)}(1+x)^{2(n-k)}=\sum_{k=0}^{n-1}C_{k+1}{n-1 \choose
k}x^k(1+x)^{k}.
\end{equation}
\end{theo}

\pf Given a plane tree $T$ with $n$ edges, we assign the weights
to the edges of $T$ by the following rule: All the terminal edges
except  the critical edge are assigned the weight $x^2$, all the
non-terminal edges are given the weight  $(1+x)^2$, and the
critical edge is assigned the weight $1$. Then the bijection
$\Phi$  transform $T$ into a  $2$-Motzkin path in which all the
down steps and wavy level steps have the weight $x^2$ and all the
up steps and straight level steps  have the weight $(1+x)^2$.

By the above weight assignment, the left hand side of
(\ref{identity2}) is then the sum of the weights of all plane
trees with $n$ edges. We now proceed to show  that the right hand
side of (\ref{identity2}) is the sum of weights of all $2$-Motzkin
paths of length $n-1$. Consider a $2$-Motzkin paths that has $k$
up steps and $k$ down steps. Since the up steps have weight $x^2$
and the down steps have weight $(1+x)^2$, it makes no difference
with respect to the sum of weights if one changes the weights of
both up steps and down steps to $x(1+x)$. In other words, such an
operation on the change of weights  is a weight-preserving
bijection on the set of $2$-Motzkin paths.

Note that for any weight assignment, we may  transform the sum of
weights of all $2$-Motzkin paths of length $n-1$ to the sum of
weights of all Motzkin paths of the same length by the following
weight assignment: the up steps and down steps carry the same
weight, and the horizontal steps in the Motzkin paths carry the
weight as the sum of the weights of a straight level step and a
wavy level step in the $2$-Motzkin path. Therefore, for our weight
assignment the sum of weights of $2$-Motzkin paths of length $n-1$
equals the sum of Motzkin paths of length $n-1$ given the
following weight assignment: up steps and down steps are given the
weight $x(1+x)$, and the horizontal steps are given the weight
\[ x^2+(1+x)^2= 1 +  x(1+x) + x(1+x).\]
So we have transformed the sum of weights of $2$-Motzkin paths of
length $n-1$ to the sum of weights of all Motzkin paths of length
$n-1$ in which all the up steps, down steps have weights $x(1+x)$,
and the horizontal steps can be regarded as either a straight
level step with weight $x(1+x)$, or a wavy level step with weight
$x(1+x)$ or a special dotted step with weight $1$.

We now get the desired sum as on the right hand side of
(\ref{identity2}) by considering the distribution of the special
dotted steps, because the remaining steps (up, down, straight
level, wavy level) all have the weight $x(1+x)$ and they  form a
$2$-Motzkin path. \qed

Setting $x={1/4}$ in  (\ref{identity1}) we obtain (\ref{q1}).


 \vskip 5mm

\noindent{\bf Acknowledgments.} This work was done under the
auspices of the National Science Foundation, the Ministry of
Education, and the Ministry of Science and Technology of China.

\small

\end{document}